\documentclass[10pt,draft,twoside,leqno,a4paper]{article}
\usepackage{amsmath,amsthm,amsfonts,amssymb,verbatim}
\newtheorem{Theorem}{Theorem}
\newtheorem{Lemma}{Lemma}

\newtheorem{Proposition}{Proposition}

\newcommand{\dash}{\mathchoice
    {\mkern.70mu\raise.50ex\hbox{\vrule height.1ex width.40em depth0pt}}
    {\mkern.40mu\raise.48ex\hbox{\vrule height.1ex width.30em depth0pt}}
    {\mkern.33mu\raise.30ex\hbox{\vrule height.1ex width.25em depth0pt}}
    {\mkern.10mu\raise.20ex\hbox{\vrule height.1ex width.20em depth0pt}}
    }

\newcommand{\R}{\mathbb R}

\newcommand{\eps}{\varepsilon}
\newcommand{\de}{\partial}
\renewcommand{\d}{{\mathrm d}}

\newcommand{\bra}{\langle}
\newcommand{\ket}{\rangle}

\newcommand{\onabla}{\overline\nabla}
\newcommand{\Cab}{C(a,b)}

\newcommand{\gx}{g_{x}}
\newcommand{\Sr}{{S_\rho(x)}}

\newcommand{\distx}{d_x}
\renewcommand{\d}{\,{\mathrm d}}
\renewcommand{\onabla}{\overline\nabla}
\begin{document}
\title{On the best H\"older exponent for
two dimensional elliptic equations
in divergence form}
\author{Tonia Ricciardi\thanks{
Supported in part by 
the MIUR National Project {\em Variational Methods and
Nonlinear Differential Equations}.}\\
{\small Dipartimento di Matematica e Applicazioni
``R.~Caccioppoli"}\\
{\small Universit\`a di Napoli Federico II}\\
{\small Via Cintia, 80126 Naples, Italy}\\
{\small fax: +39 081 675665}\\
{\small\tt{ tonia.ricciardi@unina.it}}\\
}
\date{}
\maketitle
\begin{abstract}
We obtain an estimate for the H\"older continuity
exponent for
weak solutions to the following elliptic equation 
in divergence form:
\[
\mathrm{div}(A(x)\nabla u)=0
\qquad\mathrm{in\ }\Omega,
\]
where $\Omega$ is a bounded open subset of $\R^2$
and, for every $x\in\Omega$,
$A(x)$ is a matrix with bounded measurable coefficients.
Such an estimate ``interpolates" between the 
well-known estimate of Piccinini and Spagnolo
in the isotropic case $A(x)=a(x)I$, where $a$
is a bounded measurable function,
and our previous result in the unit determinant case $\det A(x)\equiv1$.
Furthermore, we show that our estimate is sharp. Indeed, for every $\tau\in[0,1]$ we construct
coefficient matrices $A_\tau$ such that $A_0$ is isotropic and
$A_1$ has unit determinant, and such that our estimate for $A_\tau$
reduces to an equality, for every $\tau\in[0,1]$.
\end{abstract}
\begin{description}
\item {\textsc{Key Words:}} linear elliptic equation, measurable coefficients,
H\"older regularity
\item {\textsc{MSC 2000 Subject Classification:}} 35J60
\end{description}
\section{Introduction and main results}
Let $\Omega$ be a bounded open subset of $\R^2$
and let $u\in H_{\mathrm{loc}}^1(\Omega)$ be a weak solution to the
elliptic, divergence-form equation with measurable coefficients:
\begin{equation}
\label{elliptic}
\mathrm{div}(A(x)\nabla u)=0
\qquad\qquad\mathrm{in}\ \Omega,
\end{equation}
where $A(x)$, $x\in\Omega$,
is a $2\times2$ matrix satisfying the uniform ellipticity condition
\begin{equation} 
\label{ellipticity}
\lambda|\xi|^2\le\bra\xi,A(x)\xi\ket\le\Lambda|\xi|^2,
\end{equation}
for every $x\in\Omega$, for all $\xi\in\R^2$ and for some constants $0<\lambda\le\Lambda$.
By classical results of De Giorgi~\cite{DG} and Nash~\cite{Na},
it is well-known that $u$ is locally H\"older continuous in $\Omega$.
Namely, there exists a constant $\alpha\in(0,1)$ such that
for every $K\Subset\Omega$ there exists $C(K)>0$ such that:
\begin{equation}
\label{holderdef}
\frac{|u(x)-u(y)|}{|x-y|^\alpha}\le C(K)
\qquad\qquad\forall x,y\in K.
\end{equation}
The sharp estimate of $\alpha$ in terms of the ellipticity constant
$L=\Lambda/\lambda$ was obtained by Piccinini and Spagnolo~\cite{PS},
who showed that $\alpha\ge L^{-1/2}$. 
Under additional assumptions on $A$, this estimate may be improved.
For example, if
$A$ is isotropic, namely if $A(x)=a(x)I$ for some measurable function $a$
satisfying $1\le a\le L$, it was shown by Piccinini and Spagnolo~\cite{PS} that 
\begin{equation}
\label{PSalpha}
\alpha\ge\frac{4}{\pi}\arctan L^{-1/2}.
\end{equation}
On the other hand, we showed in \cite{sharpholder} that if $A$ has unit determinant, namely if
$\det A(x)=1$ for all $x\in\Omega$, then
\begin{equation}
\label{sharpholderalpha}
\alpha\ge\left(\sup_{S_\rho(x)\subset\Omega}\int_{S_\rho(x)}\bra n,A\,n\ket
\right)^{-1},
\end{equation}
where $S_\rho(x)$ is the circle of radius $\rho$ centered at $x$
and $n$ is the outward unit normal.
See Iwaniec and Sbordone~\cite{IS} for the relevance of 
the unit determinant case in the context of quasiconformal mappings.
Our aim in this note is to obtain an estimate for $\alpha$ 
in the case of general coefficient matrices $A$
satisfying the ellipticity condition \eqref{ellipticity},
which ``unifies" the estimates \eqref{PSalpha}--\eqref{sharpholderalpha}.
We shall obtain a formula which, despite of its complicated form, is indeed attained 
on a family of coefficient matrices $A_\tau$, $\tau\in[0,1]$,
such that
$A_0$ is isotropic and $A_1$ has unit determinant. 
In fact, our main effort in this note is to construct $A_\tau$.
More precisely, for every $A$ satisfying \eqref{elliptic}, let
\[
\alpha(A)=\sup\left\{\alpha\in(0,1)\ :\ 
\begin{matrix}
\mathrm{property\ \eqref{holderdef}\ holds\ for\ every}\ \\
\mathrm{solution}\ u\in H_{\mathrm{loc}}^1(\Omega)\ \mathrm{to\ \eqref{elliptic}}
\end{matrix}\right\}.
\]
We prove the following results.
\begin{Theorem}[Estimate]
\label{thm:estimate}
Suppose $A$ satisfies \eqref{ellipticity}.
Then, $\alpha(A)\ge\beta(A)$, where
\begin{equation}
\label{baralpha}
\beta(A)=\left(
\sup_{S_\rho(x)\subset\Omega}\frac{\frac{1}{|S_\rho(x)|}\int_{S_\rho(x)}
\frac{\bra n,A\,n\ket}{\sqrt{\det A}}}
{\frac{4}{\pi}\arctan\left(\frac{\inf_{S_\rho(x)}\det A}
{\sup_{S_\rho(x)}\det A}\right)^{1/4}}\right)^{-1}.
\end{equation}
\end{Theorem}
As already mentioned, Theorem~\ref{thm:estimate} is sharp, in the following sense.
\begin{Theorem}[Sharpness]
\label{thm:sharp}
For every $\tau\in[0,1]$ and for every $x\neq0$, let
$A_\tau=JK_\tau J^*$,
where
\[
K(x)=\begin{cases}
\mathrm{Id_{\R^2}},&\mathrm{if\ }\arg x\in[0,\frac{\pi}{1+M^{-\tau}})
\cup[\pi,\pi+\frac{\pi}{1+M^{-\tau}})\\
\mathrm{diag}(M,M^{1-2\tau}),&\mathrm{otherwise}
\end{cases},
\]
for some $M>1$ and
\begin{equation*}
J(\arg x)=\left(
\begin{matrix}
\cos(\arg x)&-\sin(\arg x)\\
\sin(\arg x)&\cos(\arg x)
\end{matrix}\right).
\end{equation*}
There exists $m_0>1$ such that the equality
\[
\alpha(A_\tau)=\beta(A_\tau)
\]
holds for all $M\in(1,m_0^{1/\tau})$ if $\tau>0$,
and with no restriction on $M$ if $\tau=0$.
\end{Theorem}
\subparagraph*{Notation}
Throughout this note, for all $x\in\R^2$ and for all $\rho>0$,
$B_\rho(x)$ denotes the ball of radius $\rho$ centered at $x$
and $S_\rho(x)=\de B_\rho(x)$. For every curve $\gamma$
we denote by $|\gamma|$ the length of $\gamma$. 
For every measurable function $f$
we denote by $\inf f$ and $\sup f$
the essential lower bound and the essential upper bound of $f$,
respectively.
\section{Proof of Theorem~\ref{thm:estimate}}
The proof of Theorem~\ref{thm:estimate} relies on 
an argument of Piccinini and Spagnolo~\cite{PS},
together with a weighted Wirtinger inequality obtained in \cite{wirtinger}.
Let
\[
\mathcal B=\{a\in L^\infty(\R)\ :\ a\ \mathrm{is}\ 2\pi\mathrm{-periodic\ and\ }\inf a>0\}
\]
and for every $L>1$, let
Let $\Cab>0$ denote the best constant in the
following weighted Wirtinger type
inequality:
\begin{align}
\label{wirtinger}
\int_0^{2\pi}a\,w^2\le\Cab\int_0^{2\pi}bw'^2,
\end{align}
where $w\in H_{\mathrm{loc}}^1(\R)$ is $2\pi$-periodic and satisfies
the constraint
\begin{equation}
\label{constraint}
\int_0^{2\pi}aw=0,
\end{equation}
and $a,b\in\mathcal B$. 
\begin{Lemma}[\cite{wirtinger}]
\label{lem:wirtinger}
Suppose $a,b\in\mathcal B$. Then,
\begin{equation}
\label{prelimCab}
\Cab\le
\left(
\frac{\frac{1}{2\pi}\int_0^{2\pi}\sqrt{ab^{-1}}}
{\frac{4}{\pi}\arctan\left({\frac{\inf{ab}}
{\sup{ab}}}\right)^{1/4}}\right)^2.
\end{equation}
\end{Lemma}
We note that Lemma~\ref{lem:wirtinger} reduces to the sharp 
Wirtinger inequality of Piccinini and Spagnolo~\cite{PS} when $a=b$.
Estimate \eqref{prelimCab} has been recently extended in \cite{Gi} to the case
$a,b^{-1}, \sqrt{ab^{-1}}\in L^1$ and $0<\inf(ab)\le\sup(ab)<+\infty$.
\par
In order to proceed, for every fixed $x\in\Omega$ and for every $\rho\in(0,d_x)$,
we denote by
$y=x+\rho e^{it}$ the polar coordinate transformation centered at $x$.
We denote 
\[
\onabla u=\left(u_\rho,\frac{u_t}{\rho}\right) 
\]
and for every $t\in\R$ we denote
\begin{equation}
\label{Jdef}
J(t)=\left(
\begin{matrix}
\cos t&-\sin t\\
\sin t&\cos t
\end{matrix}\right).
\end{equation}
Then,
\begin{equation}
\label{onabla}
\nabla u=J(\theta)\onabla u.
\end{equation}
\begin{Lemma}
\label{lem:PSargument}
For every matrix $A$ satisfying \eqref{ellipticity}, 
and for every $x\in\Omega$, $0<\rho<d_x$,
let
\[
C_A(x,\rho)=C\left(\bra e^{it},A(x+\rho e^{it})e^{it}\ket,
\frac{\det A(x+\rho e^{it})}{\bra e^{it},A(x+\rho e^{it})e^{it}\ket}\right)
\]
denote the best constant in \eqref{wirtinger}--\eqref{constraint} with
\begin{align*}
&a(t)=\bra e^{it},A(x+\rho e^{it})e^{it}\ket,
&&b(t)=\frac{\det A(x+\rho e^{it})}{\bra e^{it},A(x+\rho e^{it})e^{it}\ket}.
\end{align*}
Then, $\alpha(A)\ge\beta_0(A)$, where
\[
\beta_0(A)=\left(
\sup_{x\in\Omega,0<\rho<\distx}C_A(x,\rho)^{1/2}
\right)^{-1}.
\]
\end{Lemma}
\begin{proof}
We show that for every $u\in H_{\mathrm{loc}}^1(\Omega)$
solution to \eqref{elliptic} there holds:
\begin{equation}
\label{morrey}
\sup_{0<\rho<\distx}\rho^{-2\beta_0(A)}\int_{B_\rho(x)}
\bra\nabla u,A\nabla u\ket
<+\infty,
\end{equation}
for every $x\in\Omega$. 
Once estimate \eqref{morrey} is established, the statement
follows by the well-known regularity results of Morrey~\cite{Mo}.
In order to derive \eqref{morrey}, we exploit some ideas in \cite{PS}.
For every $x\in\Omega$ and for every $0<\rho<\distx$, we set
\[
\gx(\rho)=\int_{B_\rho(x)}
\bra\nabla u,A\nabla u\ket.
\]
We denote by $P=(p_{ij})$ the matrix defined by
\[
P(x+\rho e^{it})=J^*(t)A(x+\rho e^{it})J(t).
\]
Note that $p_{11}(x+\rho e^{it})=\bra e^{it},A(x+\rho e^{it})e^{it}\ket$
and $\det P=\det A$.
By the divergence theorem and \eqref{elliptic}, we have
\[
\gx(\rho)=\int_\Sr(u-\mu)\bra n,A\nabla u\ket
=\int_\Sr(u-\mu)\bra\underline e_1,P\onabla u\ket,
\]
where 
$n$ is the outward normal to $\Sr$, $\underline e_1=(1,0)$
and $\mu$ is any constant.
In view of H\"older's inequality, we may write
\begin{align*}
\gx(\rho)\le\left(
\int_\Sr p_{11}(u-\mu)^2\right)^{1/2}
\left(\int_\Sr\frac{\bra\underline e_1,P\onabla u\ket^2}{p_{11}}\right)^{1/2}.
\end{align*}
By inequality \eqref{wirtinger} with 
$a(t)=p_{11}(x+\rho e^{it})$, $b(t)=\det A/p_{11}(x+\rho e^{it})$ and
\begin{align*}
&w(t)=u(x+\rho e^{it})-\mu, 
&\mu=\frac{1}{2\pi}\int_0^{2\pi}p_{11}(x+\rho e^{it})u(x+\rho e^{it})\d\theta,
\end{align*}
we derive 
\[
\int_{\Sr}p_{11}(u-\mu)^2\le C_A(x,\rho)\int_{\Sr}\frac{\det A}{p_{11}}u_t^2.
\]
Therefore,
\begin{align*}
\gx(\rho)\le C_A^{1/2}(x,\rho)\left(
\int_\Sr \frac{\det A}{p_{11}}u_t^2\right)^{1/2}
\left(\int_\Sr\frac{\bra\underline e_1,P\onabla u\ket^2}{p_{11}}\right)^{1/2}.
\end{align*}
At this point, we observe that any $2\times2$ symmetric matrix $B=(b_{ij})$
such that $b_{11}\neq0$
satisfies the following identity:
\begin{equation}
\label{matrixidentity}
\bra\xi,B\xi\ket
=\frac{\bra\xi,B e_1\ket^2}{b_{11}}+\frac{\det B}{b_{11}}\bra\xi,e_2\ket^2,
\end{equation}
for any $\xi\in\R^2$.
Recalling that $u_\theta/\rho=(\onabla u)_{22}$, 
in view of the elementary inequality $\sqrt{ab}\le(a+b)/2$ and
of identity \eqref{matrixidentity}
with $B=P(x+\rho e^{i\theta})$ and $\xi=\onabla u$,
we obtain:
\begin{align*}
\gx(\rho)\le&\rho\,C_A^{1/2}(x,\rho)\left(
\int_\Sr \frac{\det A}{p_{11}}
\left(\frac{u_t}{\rho}\right)^2\right)^{1/2}
\left(\int_\Sr\frac{\bra e_1,P\onabla u\ket^2}{p_{11}}\right)^{1/2}\\
=&\rho\,C_A^{1/2}(x,\rho)\left(
\int_\Sr \frac{\det A}{p_{11}}
\left(\onabla u\right)_{22}^2\right)^{1/2}
\left(\int_\Sr\frac{\bra e_1,P\onabla u\ket^2}{p_{11}}\right)^{1/2}\\
\le&\frac{\rho}{2}\,C_A^{1/2}(x,\rho)
\int_\Sr\left(\frac{\det A}{p_{11}}
\left(\onabla u\right)_{22}^2
+\frac{\bra e_1,P\onabla u\ket^2}{p_{11}}\right)\\
=&\frac{\rho}{2}\,C_A^{1/2}(x,\rho)
\int_\Sr\bra\onabla u,P\onabla u\ket
=\frac{\rho\,C_A^{1/2}(x,\rho)}{2}
\int_\Sr\bra\nabla u,A\nabla u\ket.
\end{align*}
Recalling the definition of $\gx$, we derive from the above inequality that:
\[
\gx(\rho)\le\frac{\rho}{2}\,C_A^{1/2}(x,\rho)\,\gx'(\rho),
\]
for almost every $0<\rho<d_x$.
In particular, for every $\rho\in(0,d_x)$ we have:
\[
\gx(\rho)\le\frac{\rho}{2}\sup_{x\in\Omega,0<\rho<\distx}C_A^{1/2}(x,\rho)\,\gx'(\rho)
=\frac{\rho}{2\beta_0(A)}\,\gx'(\rho)
\qquad\mathrm{in\ }(0,\distx).
\]
The above implies that the function
$\rho^{-2\beta_0(A)}\gx(\rho)$ in non-decreasing, and therefore bounded, 
in $(0,\distx)$.
\end{proof}
\begin{proof}[Proof of Theorem~\ref{thm:estimate}]
In view of Lemma~\ref{lem:wirtinger}--(i)
with
\[
a(t)=p_{11}(x+\rho e^{it}),
\qquad
b(t)=\frac{\det A}{p_{11}}(x+\rho e^{it}),
\]
we have
\[
C_A(x,\rho)\le
\left(\frac{\frac{1}{2\pi}\int_0^{2\pi}\frac{p_{11}}{\sqrt{\det A}}(x+\rho e^{it})\d t}
{\frac{4}{\pi}\arctan\left(\frac{\inf_{t\in(0,2\pi)}\det A(x+\rho e^{it})}
{\sup_{t\in(0,2\pi)}\det A(x+\rho e^{it})}\right)^{1/4}}\right)^2.
\]
Now the asserted estimate follows by Lemma~\ref{lem:PSargument}.
The sharpness will follow by Proposition~\ref{prop:sharp} in the
next section.
\end{proof}
\section{Proof of Theorem~\ref{thm:sharp}}
We define
\begin{align*}
&\mu=\frac{4}{\pi}\arctan M^{-(1-\tau)/2},&&c=\frac{2}{1+M^{-\tau}}.
\end{align*}
We prove
\begin{Proposition}
\label{prop:sharp}
For every $\tau\in[0,1]$, let $A_\tau$ be the coefficient
matrix defined in Theorem~\ref{thm:sharp}.
There exists $m_0>1$ such that 
there holds 
\begin{align*}
\beta(A_\tau)=\frac{\mu}{c}=
\frac{2}{\pi}
(1+M^{-\tau})\arctan M^{-(1-\tau)/2}
\end{align*}
for all $M\in(1,m_0^{1/\tau})$ if $\tau>0$, and with no restriction on $M$
if $\tau=0$.
Furthermore, let $u_\tau\in H^1(B)$ be defined in polar coordinates by
\[
u_\tau(\rho,\theta)=\rho^{\mu/c}w_\tau(\theta),
\]
where 
\[
w_\tau(\theta)=\begin{cases}
\sin[\mu(c^{-1}\theta-\pi/4)],&\mathrm{if\ }\theta\in[0,c\pi/2)\\
M^{-(1-\tau)/2}\cos[\mu\left(c^{-1}M^{\tau}(\theta-c\pi/2)-\pi/4\right)],
&\mathrm{if\ }\theta\in[c\pi/2,\pi)
\end{cases}
\]
and $w_\tau(\pi+\theta)=-w_\tau(\theta)$ for all $\theta\in[0,\pi)$.
Then, $u_\tau$
is a weak solution to the elliptic equation \eqref{elliptic} 
with $A=A_\tau$, 
and its H\"older exponent is $\mu/c$.
In particular, $\alpha(A_\tau)=\beta(A_\tau)$.
\end{Proposition}
In order to prove Proposition~\ref{prop:sharp},
we begin by proving some lemmas.
\begin{Lemma}
\label{lem:polar}
For every $x\neq0$ let $\theta=\arg x$
and let
\[
A(x)=A(\theta)=J(\theta)K(\theta)J^*(\theta),
\]
where  
\[
K(\theta)=\left(
\begin{matrix}
k_1(\theta)&0\\
0&k_2(\theta)
\end{matrix}
\right)
\]
for some positive and bounded, $2\pi$-periodic functions $k_1,k_2$.
Then, in polar coordinates, equation \eqref{elliptic} takes the form:
\begin{equation}
\label{polar}
\begin{cases}
\left(\rho k_1u_\rho\right)_\rho+\left(\frac{k_2}{\rho}u_\theta\right)_\theta=0
\qquad\mathrm{in\ }(0,+\infty)\times\R\\
u\ 2\pi\mathrm{-periodic\ in\ }\theta
\end{cases}.
\end{equation}
If $u\not\equiv0$ is of the separation of variables form $u(\rho,\theta)=R(\rho)\Theta(\theta)$, then
$u$ satisfies \eqref{polar} if and only if
$R(\rho)=\rho^\gamma$ for some constant $\gamma>0$
and $\Theta$ is a $2\pi$-periodic weak solution to
the equation:
\begin{align}
\label{sep2}
-(k_2\Theta')'=\gamma^2 k_1\Theta
\qquad\mathrm{in\ }\R.
\end{align}
\end{Lemma}
\begin{proof}
By definition, $u$ satisfies
\[
\int_{\R^2}\bra\nabla u,A\nabla v\ket=0
\qquad\forall v\in C_c^\infty(\R^2).
\]
In polar coordinates centered at 0, recalling that $\onabla u=(u_\rho,u_\theta/\rho)$ , we have:
\begin{align*}
0=\int_B\bra\nabla u,A\nabla v\ket
=&\int_{(0,+\infty)\times(0,2\pi)}\bra\onabla u,K(\theta)\onabla v\ket\rho\d\rho\d\theta\\
=&\int_{(0,+\infty)\times(0,2\pi)}\left(\rho k_1u_\rho v_\rho
+\frac{k_2}{\rho}u_\theta v_\theta\right)\d\rho\d\theta,
\end{align*}
for every $v\in C_c^\infty(B)$.
Integration by parts yields \eqref{polar}.
Now suppose that $u(\rho,\theta)=R(\rho)\Theta(\theta)$. In view of Nikodym's theorem,
$R$ and $\Theta$ are absolutely continuous on $(0,+\infty)$ and $\R$, respectively.
Choosing $v$ of the form $v(\rho,\theta)=\varphi(\rho)\psi(\theta)$
with $\varphi\in C_c^\infty(0,+\infty)$ and $\psi\in C_c^\infty(0,2\pi)$,
we derive from the above that
\begin{align*}
\int_0^{+\infty}\rho R'\varphi'\d\rho\int_0^{2\pi}k_1\Theta\psi\d\theta
+\int_0^{+\infty}\frac{R}{\rho}\varphi\d\rho\int_0^{2\pi}k_2\Theta'\psi'\d\theta=0.
\end{align*}
Since $\varphi,\psi$ are arbitrary, we conclude that
\begin{align*}
\frac{\int_0^{+\infty}(\rho R')'\varphi\d\rho}{\int_0^{+\infty}R\rho^{-1}\varphi\d\rho}
=-\frac{\int_0^{2\pi}(k_2\Theta')'\psi\d\theta}{\int_0^{2\pi}k_1\Theta\psi\d\theta}=\tau,
\end{align*}
for some constant $\tau\in\R$.
It follows that
\[
\int_0^{+\infty}(\rho R')'\varphi\d\rho=\tau\int_0^{+\infty}\frac{R}{\rho}\varphi\d\rho
\qquad\forall\varphi\in C_c^\infty(0,+\infty)
\]
and
\[
\int_0^{2\pi}(k_2\Theta')'\psi\d\theta=-\tau\int_0^{2\pi}k_1\Theta\psi\d\theta
\qquad\forall\psi\in C_c^\infty(0,2\pi),
\]
and therefore \eqref{sep2} is established.
By regularity, $R$ is smooth in $(0,+\infty)$, it
satisfies $(\rho R')'=\tau R\rho^{-1}$ in $(0,+\infty)$ and
is bounded at the origin. 
Therefore, $R(\rho)=\rho^\gamma$
with $\gamma^2=\tau>0$.
\end{proof}
\begin{Lemma}
\label{lem:expansion}
Suppose $A$ satisfies the assumptions of Lemma~\ref{lem:polar}.
Then, for all $x\in\R^2$, $\rho>0$ and $t\in\R$ such that
$x+\rho e^{it}\neq0$, we have
\begin{align}
\label{bracket}
\frac{\bra e^{it},A(x+\rho e^{it})e^{it}\ket}
{\sqrt{\det{A(x+\rho e^{it})}}}
=\sqrt{\frac{k_1(\theta(t))}{k_2(\theta(t))}}\cos^2\left(\theta(t)-t\right)
+\sqrt{\frac{k_2(\theta(t))}{k_1(\theta(t))}}\sin^2\left(\theta(t)-t\right),
\end{align}
where
\[
\theta(t)=\arg(x+\rho e^{it}).
\]
\end{Lemma}
\begin{proof}
Using the fact that $J^*(\theta)e^{it}=e^{i(t-\theta)}$ for all $t,\theta\in\R$, we have:
\begin{align*}
\bra e^{it},A(x+\rho e^{it})e^{it}\ket=&\bra J^*(\theta(t))e^{it},K(\theta(t))J^*(\theta(t))e^{it}\ket\\
=&k_1(\theta(t))\cos^2(t-\theta(t))+k_2(\theta(t))\sin^2(t-\theta(t)).
\end{align*}
Now \eqref{bracket} follows easily.
\end{proof}
We shall need the following property from Euclidean geometry.
As we have not found a proof in the literature, we include one here.
\begin{Lemma}
\label{lem:euclid}
Let $\mathcal C$ be a (two-sided) cone with vertex at the origin and let 
$x\in\R^2$ be such that $|x|<1$.  Then
\begin{equation}
\label{euclid}
|\mathcal C\cap S_1(x)|=|\mathcal C\cap S_1(0)|.
\end{equation}
\end{Lemma}
\begin{proof}
We denote by $A,B,C,D$ the intersection points of $\mathcal C$ with $S_1(x)$
taken in, say, counterclockwise order.
We have to show that
$\angle AxB+\angle CxD=\angle AOB+\angle COD=2\angle AOB$.
We set $\alpha=\angle AxB$, $\beta=\angle CxD$, $\eps=\angle xAC=\angle xCA$,
$\delta=\angle xBD=\angle xDB$, $\eta=\angle ABx=\angle BAx$, $\theta=\angle xCD=\angle xDC$,
$\varphi=\angle AOB=\angle COD$. Then, summing the angles of the triangles
$AxB$, $CxD$, $AOB$, $COD$, respectively, we obtain:
\begin{align*}
&\alpha+2\eta=\pi,&&\eta+\delta+\eta-\eps+\varphi=\pi\\
&\beta+2\theta=\pi,
&&\theta+\eps+\theta-\delta+\varphi=\pi.
\end{align*}
Summation of these equations yields
$\alpha+\beta=2\pi-2(\eta+\theta)$
and $2(\eta+\theta)=2\pi-2\varphi$, from which we derive
the desired equality $\alpha+\beta=2\varphi$.
\end{proof}
For every $x\in\R^2$ and for every $\rho>0$ we define
\[
f(x,\rho)=\frac{1}{|S_\rho(x)|}\int_{S_\rho(x)}\frac{\bra n,A_\tau\,n\ket}
{\sqrt{\det{A_\tau}}}
=\frac{1}{2\pi}\int_0^{2\pi}\frac{\bra e^{it},A_\tau(x+\rho e^{it})e^{it}\ket}
{\sqrt{\det{A_\tau(x+\rho e^{it})}}}\d t,
\]
where $A_\tau$ is the matrix defined in Theorem~\ref{thm:sharp}. We note that
\begin{equation}
\label{homog}
f(x,\rho)=f\left(\frac{x}{\rho},1\right).
\end{equation}
We prove the following.
\begin{Lemma}
\label{lem:fest}
There exists $m_0>1$ such that for all $x\in\R^2$ and for all $\rho>0$
there holds
\[
f(x,\rho)\le f(0,1)=c=\frac{2}{1+M^{-\tau}},
\]
for all $M\in(1,m_0^{1/\tau})$ if $\tau>0$, and with no restriction on $M$
if $\tau=0$.
\end{Lemma}
\begin{proof}
Throughout this proof, we let
\[
m:=M^{\tau}
\]
and
\[
\mathcal C:=\left\{x\in\R^2\setminus\{0\}:\ \arg x\in[\frac{\pi}{2}c,\pi)\cup[\pi+\frac{\pi}{2}c,2\pi)\right\}.
\]
Then
\[
K_\tau(x)=\begin{cases}
\mathrm{diag}(M,M^{1-2\tau}),&\mathrm{if\ }x\in\mathcal C\\
\mathrm{Id_{\R^2}},&\mathrm{otherwise}
\end{cases}.
\]
In view of Lemma~\ref{lem:expansion}, it follows that
\begin{align*}
\label{bracket}
&\frac{\bra e^{it},A_\tau(x+\rho e^{it})e^{it}\ket}
{\sqrt{\det{A_\tau(x+\rho e^{it})}}}=\\
&\qquad=\begin{cases}
m\cos^2\left(\theta(t)-t\right)+m^{-1}\sin^2\left(\theta(t)-t\right),
&\mathrm{if}\ x+\rho e^{it}\in\mathcal C\\
1,
&\mathrm{otherwise}
\end{cases}
\end{align*}
and in view of \eqref{homog}, we may assume $\rho=1$.
\par
Case (i): $|x|<1$.
We estimate:
\begin{align*}
2\pi f(x,1)=&2\pi-|\mathcal C\cap S_1(x)|\\
&+\int_{\mathcal C\cap S_1(x)}
\left\{m\cos^2\left(\theta(t)-t\right)
+\frac{1}{m}\sin^2\left(\theta(t)-t\right)\right\}\d t\\
=&2\pi-|\mathcal C\cap S_1(x)|
+\int_{\mathcal C\cap S_1(x)}
\left\{\left(m-\frac{1}{m}\right)\cos^2\left(\theta(t)-t\right)+\frac{1}{m}\right\}\d t\\
\le&2\pi-|\mathcal C\cap S_1(x)|+m|\mathcal C\cap S_1(x)|
=2\pi+(m-1)|\mathcal C\cap S_1(x)|.
\end{align*}
From Lemma~\ref{lem:euclid} we derive $|\mathcal C\cap S_1(x)|=|\mathcal C\cap S_1(0)|$
and therefore we obtain the estimate
\begin{align*}
2\pi f(x,1)\le2\pi+(m-1)|\mathcal C\cap S_1(0)|.
\end{align*}
On the other hand, when $x=0$ we have $\theta(t)-t\equiv0$, and 
consequently
\begin{align*}
\frac{\bra e^{it},A_\tau(e^{it})e^{it}\ket}{\sqrt{\det A_\tau(e^{it})}}
=\begin{cases}
m,&\mathrm{if\ }e^{it}\in \mathcal C\\
1&\mathrm{otherwise}
\end{cases}.
\end{align*}
It follows that 
$2\pi f(0,1)=2\pi+(m-1)|\mathcal C\cap S_1(0)|$,
and recalling the definition of $\mathcal C$ and $c$,
we obtain
$f(0,1)=c=2/(1+m^{-1})$.
Hence, the desired estimate 
$f(x,1)\le f(0,1)$ follows in the case $|x|<1$
with no restriction on $M$.
\par
Case (ii): $|x|\ge1$.
We set $h(t)=\cos^2(\theta(t)-t)$.
By elementary geometrical arguments, for every $0\le k\le 1$ 
we have
\begin{equation}
\label{arcsin}
\left|\{h(t)\le k\}\right|=4\arcsin\sqrt k.
\end{equation}
Since $mh(t)+m^{-1}(1-h(t))\ge1$ if and only if $h(t)\ge(m+1)^{-1}$,
we estimate
\[
2\pi f(x,1)
\le\left|\left\{h(t)\le\frac{1}{m+1}\right\}\right|
+\int_{\{h(t)\ge(m+1)^{-1}\}}\Big\{mh(t)+\frac{1}{m}(1-h(t))\Big\}\d t.
\]
By virtue of \eqref{arcsin}, we derive
\begin{equation}
\label{ffirst}
2\pi f(x,1)
\le4\arcsin\sqrt{\frac{1}{m+1}}
+\int_{\{h(t)\ge(m+1)^{-1}\}}\Big\{mh(t)+\frac{1}{m}(1-h(t))\Big\}\d t.
\end{equation}
Similarly, let $\eps>0$ and note that $1+\eps(m-1)/m\le m$. We have that
$mh(t)+m^{-1}(1-h(t))\ge1+\eps(m-1)/m$ if and only if $h\ge(1+\eps)/(m+1)$.
Therefore, we estimate in turn
\begin{align*}
&\int_{\{h(t)\ge(m+1)^{-1}\}}\Big\{mh(t)+\frac{1}{m}(1-h(t))\Big\}\d t\\
&\le\left(1+\eps\frac{m-1}{m}\right)
\left|\left\{\frac{1}{m+1}\le h\le\frac{1+\eps}{m+1}\right\}\right|
+m\left(2\pi-\left|\left\{h\le\frac{1+\eps}{m+1}\right\}\right|\right)\\
&=4\left(1+\eps\frac{m-1}{m}\right)\left(\arcsin\sqrt{\frac{1+\eps}{m+1}}
-\arcsin\sqrt{\frac{1}{m+1}}\right)\\
&\qquad\qquad+m\left(2\pi-4\arcsin\sqrt{\frac{1+\eps}{m+1}}\right)\\
&=2\pi m-4\left(1+\eps\frac{m-1}{m}\right)\arcsin\sqrt{\frac{1}{m+1}}\\
&\qquad\qquad-4(m-1)\left(1-\frac{\eps}{m}\right)\arcsin\sqrt{\frac{1+\eps}{m+1}}.
\end{align*}
Hence, 
\begin{align*}
f(x,1)\le
2\pi m-4\eps\frac{m-1}{m}\arcsin\sqrt{\frac{1}{m+1}}
-4(m-1)\left(1-\frac{\eps}{m}\right)\arcsin\sqrt{\frac{1+\eps}{m+1}},
\end{align*}
and it suffices to check that there exist $\eps>0$ and $m_0>1$ such that
\begin{align*}
\frac{1}{2\pi}&\left(2\pi m-4\eps\frac{m-1}{m}\arcsin\sqrt{\frac{1}{m+1}}
-4(m-1)\left(1-\frac{\eps}{m}\right)\arcsin\sqrt{\frac{1+\eps}{m+1}}
\right)\\
&\qquad\qquad\qquad\le\frac{2}{1+m^{-1}},
\end{align*}
for all $1<m\le m_0$.
Upon factorization, the above is equivalent to:
\begin{align}
\label{factor}
\frac{m(m-1)}{m+1}\le(m-1)
\left[\frac{2\eps}{\pi m}\arcsin\sqrt{\frac{1}{m+1}}-\left(1-\frac{\eps}{m}\right)
\frac{2}{\pi}\arcsin\sqrt{\frac{1+\eps}{m+1}}\right].
\end{align}
Therefore, if $\tau=0$, we have $m=1$ and \eqref{factor} holds
with no restriction on $M$.
If $\tau>0$ we have $m-1>0$ and \eqref{factor} is verified if and only if
\begin{equation}
\label{mcond}
\frac{m}{m+1}\le\frac{2}{\pi}\arcsin\sqrt{\frac{1+\eps}{m+1}}
-\frac{2\eps}{\pi m}\left(\arcsin\sqrt{\frac{1+\eps}{m+1}}
-\arcsin\sqrt{\frac{1}{m+1}}\right).
\end{equation}
Let $\delta=m-1>0$ and consider the function $\zeta$ defined by
\begin{align*}
\zeta(\eps,&\delta)=\\
=&\frac{2}{\pi}\left[\arcsin\sqrt{\frac{1+\eps}{2+\delta}}
-\frac{\eps}{1+\delta}\left(\arcsin\sqrt{\frac{1+\eps}{2+\delta}}
-\arcsin\sqrt{\frac{1}{2+\delta}}\right)\right]
-\frac{1+\delta}{2+\delta}.
\end{align*}
Then \eqref{mcond} is equivalent to $\zeta(\eps,\delta)\ge0$.
We note that
$\zeta(0,0)=\frac{2}{\pi}\arcsin\sqrt{\frac{1}{2}}-\frac{1}{2}=0$.
By Taylor's expansion, there exists $\eps_0>0$ such that
the strict inequality $\zeta(\eps_0,0)>0$ is satisfied.
Hence, by continuity, there exists $\delta_0>0$ such that
$\zeta(\eps_0,\delta)>0$ for all $\delta\in(0,\delta_0)$.
Setting $m_0=1+\delta_0$, we conclude that \eqref{mcond} is satisfied
for $\eps=\eps_0$ and for all $\delta\in(0,\delta_0)$. 
It follows that the statement of the lemma holds with $M_0=m_0^{1/\tau}$.  
\end{proof}
\begin{proof}[Proof of Proposition~\ref{prop:sharp}]
In view of Lemma~\ref{lem:fest} we have
\begin{align*}
\beta(A_\tau)=&\left(\sup_{x\in\Omega,0<\rho<\distx}
\frac{f(x,\rho)}{4\pi^{-1}\arctan M^{-(1-\tau)/2}}\right)^{-1}\\
=&\frac{2}{\pi}
(1+M^{-\tau})\arctan M^{-(1-\tau)/2}.
\end{align*}
On the other hand, by direct check we see that $w_\tau$ is a 
$2\pi$-periodic weak solution to the equation
\[
-\left(k_{\tau,2}w_\tau'\right)'
=\frac{\mu}{c}k_{\tau,1}w_\tau
\qquad\mathrm{in\ }\R,
\]
where $k_{\tau,1}$, $k_{\tau,2}$ are the diagonal entries of $K_\tau$.
It follows by Lemma~\ref{lem:polar} that $u_\tau$ satisfies \eqref{elliptic}
with $A=A_\tau$. Since $w_\tau$ is absolutely continuous,
$u_\tau$ is H\"older continuous with exponent $\beta(A_\tau)$.
\end{proof}
\begin{proof}[Proof of Theorem~\ref{thm:sharp}]
The proof is a direct consequence of Proposition~\ref{prop:sharp}.
\end{proof}


\begin{thebibliography}{99}
\bibitem{DG}
E.~De Giorgi,
Sulla differenziabilit\`a e l'analiticit\`a delle
estremali degli integrali multipli regolari, 
Mem.\ Acc.\ Sci.\ Torino Cl.\ Sci.\ Fis.\ Mat.\ Nat.\ (3)
\textbf{3} (1957), 25--43.
\bibitem{Gi}
R.~Giova,
A weighted Wirtinger inequality, 
to appear on Ricerche Mat.
\bibitem{IS}
T.~Iwaniec and C.~Sbordone,
Quasiharmonic fields, 
Ann.\ Inst.\ H.~Poincar\'e Anal.\ Non Lin\'eaire \textbf{18} No.~5 (2001), 519--572.
\bibitem{Mo}
C.B.~Morrey,
On the solutions of quasi-linear elliptic partial differential equations, 
Trans.\ Am.\ Math.\ Soc.\ \textbf{43} No.~1 (1938), 126--166.
\bibitem{Na}
J.~Nash,
Continuity of solutions of parabolic and elliptic equations, 
Amer.\ J.\  Math.\ \textbf{80} No.~5 (1958), 931--954.
\bibitem{PS} 
L.C.~Piccinini and S.~Spagnolo, 
On the H\"older continuity of solutions of second order elliptic equations
in two variables, Ann.\ Scuola Norm.\ Sup.\ Pisa \textbf{26} No.~2 (1972), 391--402.
\bibitem{sharpholder}
T.~Ricciardi,  
A sharp H\"older estimate for elliptic equations
in two variables, 
Proc. Roy. Soc. Edinburgh \textbf{135A} (2005), 165--173.
\bibitem{wirtinger} 
T.~Ricciardi,  
A sharp weighted Wirtinger inequality, 
Boll. U.M.I. (8) \textbf{8-B} (2005), 259--267.
\end{thebibliography}
\end{document}